\documentstyle[11pt]{article}

\newcommand{\qcom}{$q$-commutation relations}
\newcommand{\tq}{{\cal T}_{q}}
\newcommand{\tz}{{\cal T}_{0}}
\newcommand{\eq}{{\cal E}^{q}}
\newcommand{\eo}{{\cal E}^{0}}
\newcommand{\rqiu}{{\cal R}^{q}}
\newcommand{\ro}{{\cal R}^{0}}
\newcommand{\ltq}{{\cal L} ( \tq )}
\newcommand{\ltz}{{\cal L} ( \tz )}
\newcommand{\lt}{{\cal L} ( {\cal T} )}
\newcommand{\lvn}{{\cal L} ( \vn )}
\newcommand{\lvaq}{{\cal L} ( \vaq )}

\newcommand{\lvnq}{{\cal L} ( \vnq )}
\newcommand{\lvnqvn}{{\cal L} ( \vnq , \vn )}
\newcommand{\lvnvnq}{{\cal L} ( \vn , \vnq )}
\newcommand{\od}{{\cal O}_{d}}
\newcommand{\fiq}{\Phi_{q}}
\newcommand{\fio}{\Phi_{0}}
\newcommand{\va}{{\cal V}_\alpha}
\newcommand{\vaq}{{\cal V}_{\alpha,q}}
\newcommand{\vn}{{\cal V}_{n}}
\newcommand{\vnp}{{\cal V}_{n+1}}
\newcommand{\vz}{{\cal V}_{0}}
\newcommand{\vnq}{{\cal V}_{n,q}}
\newcommand{\vzq}{{\cal V}_{0,q}}
\newcommand{\ciclukl}{(k \rightarrow l)}

\newcommand{\cicluunuk}{(1 \rightarrow k)}
\newcommand{\siri}{1 \leq i_{1} , \ldots , i_{n} \leq d}
\newcommand{\sirxii}{\xi_{i_{1}} \otimes \cdots \otimes \xi_{i_{n}}}
\newcommand{\sirxij}{\xi_{j_{1}} \otimes \cdots \otimes \xi_{j_{n}}}
\newcommand{\sirxiidoi}{\xi_{i_{2}} \otimes \cdots \otimes \xi_{i_{n}}}
\newcommand{\sirxiihat}{\xi_{i_{1}} \otimes  \cdots \otimes
\widehat{\xi_{i_{k}}} \otimes \cdots \otimes \xi_{i_{n}}}
\newcommand{\sirxijhat}{\xi_{j_{1}} \otimes  \cdots \otimes
\widehat{\xi_{j_{k}}} \otimes \cdots \otimes \xi_{j_{n}}}
\newcommand{\xn}{\sum_{k=1}^{n} q^{k-1} (1 \rightarrow k)}
\newcommand{\prinf}{\prod_{k=1}^{\infty} \frac{1-|q|^{k}}{1+|q|^{k}}}
\newcommand{\prinfinv}{\prod_{k=1}^{\infty} \frac{1+|q|^{k}}{1-|q|^{k}}}

\begin{document}

\title{\bf On the Fock representation of the $q$-commutation relations}
\author{Ken Dykema${ }^{\sharp}$ \\ Department of Mathematics \\
University of California \\ Berkeley, CA 94720 \\
e-mail: dykema@math.berkeley.edu
\and Alexandru Nica \\ Department of Mathematics \\
University of California \\ Berkeley, CA 94720 \\
e-mail: nica@cory.berkeley.edu}
\maketitle

\vspace{1.5in}

\begin{abstract}
We consider the C*-algebra $\rqiu$ generated by the
\setcounter{page}{0}
representation of the \qcom $\  $on the twisted Fock space.
We construct a canonical unitary
$U \ (=U(q))$ from the twisted Fock space
to the usual Fock space, such that
$U \rqiu U^{*}$ contains the extended Cuntz algebra
$\ro$, for all $q \in (-1,1)$. We prove the equality
$U \rqiu U^{*} = \ro$ for $q$ satisfying:
\[
q^{2} \ < \ 1-2|q|+2|q|^{4}-2|q|^{9}+ \cdots +2(-1)^{k}|q|^{k^{2}}+ \cdots .
\]
\end{abstract}

\vspace{1.8in}

${ }^{\sharp}$Partially supported by the Fannie and John Hertz
Foundation.

\vspace{3in}
\setlength{\baselineskip}{18pt}

\section{Introduction and statement of results}

$\ \ \ \ \ $In this paper we study the C*-algebra generated by the
\setcounter{page}{1}
representation on the twisted Fock space of
the \qcom . These relations, introduced
by Greenberg \cite{G} and Bo\.{z}ejko and Speicher
\cite{BS}, provide an interpolation depending on a parameter
$q \in (-1,1)$ between the bosonic and the fermionic commutation
relations (which correspond to $q=1$ and $q=-1$, respectively).
For $q$ in [-1,1], a representation of the \qcom $\ $is of the form:
\begin{equation}
c( \xi ) c( \eta )^{*} - q c( \eta )^{*} c( \xi ) \ = \
< \xi \ | \ \eta >I, \ \  \xi , \eta \in {\cal H},
\end{equation}
where $\cal H$ is a separable Hilbert space and $c( \cdot )$ is
linear with values operators on some Hilbert space $\cal K$
(called the space of the representation). The Fock representation
of these relations is the one uniquely determined, up to unitary
equivalence, by the following condition: there exists a vacuum vector
$\Omega$ in the space of the representation, that is cyclic for
the C*-algebra generated by $\{ c( \xi ) | \xi \in {\cal H} \}$,
and such that $c( \xi ) \Omega = 0$ for every $\xi$ in $\cal H$.
The uniqueness of the Fock representation is easy to show, but
the proof of its existence is not at all trivial (see \cite{BS},
\cite{F}, \cite{Z}); the construction giving this representation
will be briefly reviewed in Section 2.1 below.

We shall consider the case when the Hilbert space $\cal H$ of
(1.1) has finite dimension $d \geq 2$; $d$ will be fixed throughout
the whole paper. Choosing an orthonormal basis of $\cal H$,
we see that the representations of the \qcom $\ $come to those
of the universal unital C*-algebra generated by $d$ elements
$a_{1}, \ldots , a_{d}$, that satisfy:
\begin{equation}
a_{i} a_{j}^{*} - qa_{j}^{*} a_{i} \ = \ \delta_{i,j} I, \ \
1 \leq i,j \leq d.
\end{equation}
We shall denote, following \cite{JWS}, this universal C*-algebra
by $\eq$ (= $\eq (d)$). Also, we shall use the following notations:
the Fock representation of the \qcom , viewed as a representation
of $\eq$, will be denoted by $\fiq$, and its space will be denoted by
$\tq$ (and called the twisted Fock space); we shall put
\begin{equation}
A_{i} \ = \ \fiq (a_{i}) \in \ltq , \ \ 1 \leq i \leq d,
\end{equation}
and we shall denote by $\rqiu$ the C*-algebra generated by
$A_{1}, \ldots ,A_{d}$ in $\ltq$. Equivalently, $\rqiu = \fiq ( \eq )$;
this C*-algebra will be our main object of investigation.

For $q=0$, we have that $a_{1}, \ldots ,a_{d}$ of (1.2)  are the adjoints of
$d$ isometries with mutually orthogonal ranges, hence $\eo$ is the
well-known extension by the compacts of the Cuntz algebra $\od$ (\cite{C});
moreover, $\fio : \eo \rightarrow \ltz$ is precisely the canonical
representation of $\eo$ on the usual Fock space
$\tz = {\cal T} = {\bigoplus}_{n=0}^{\infty} \left(
{({\bf C}^{d})}^{\otimes n} \right)$ (\cite{E}).
It is known that $\fio$ is faithful, hence (if we consider that the
``non-deformed case'' is $q=0$), both $\eq$
and $\rqiu$ can be viewed
as deformations of the extension by the compacts of $\od$.

In order to distinguish the case $q=0$, we shall write
$v_{1}, \ldots ,v_{d}$ for the $a_{1}, \ldots , a_{d}$ of (1.2)
corresponding to this case (to be very rigorous, we should have
written in (1.2)  $a_{i,q}$ instead of $a_{i}$, and then $v_{i}$
would be defined as $a_{i,0}$; however, the value of $q$ which
is considered will always be clear, and we preferred to keep the
notations simple). We denote the projection
$\sum_{i=1}^{d} v_{i}^{*} v_{i} \in \eo$ by $p$. Similarly, we
shall write $V_{1}, \ldots , V_{d}$ for the $A_{1}, \ldots , A_{d}$
of (1.3) corresponding to $q=0$, and put
$P = \sum_{i=1}^{d} V_{i}^{*} V_{i}$; then $V_{1}, \ldots , V_{d}$ are
annihilation operators on the (non-deformed) Fock space
${\cal T} = \tz$ (formula (2.2) below),
and $P$ is the projection onto the orthogonal
complement of the vacuum vector.

In \cite{JWS} it was proved that, for $|q| < \sqrt{2} -1$,
$\eq \simeq \eo$ and $\fiq$ is faithful. The proof involved
finding a positive element $\rho \in \eo$ that satisfies
the equation:
\begin{equation}
\rho^{2} \ = \ p + q \sum_{i,j =1}^{d}
(v_{i} \rho v_{j})^{*} (v_{j} \rho v_{i}),
\end{equation}
and then showing that
$a_{i} \rightarrow v_{i} \rho , \ 1 \leq i \leq d$,
gives an isomorphism between $\eq$ and $\eo$. We shall consider the
analogue on the Fock space of (1.4), i.e.:
\begin{equation}
R^{2} \ = \ P + q \sum_{i,j=1}^{d}
(V_{i} R V_{j})^{*} (V_{j} R V_{i}),
\end{equation}
($P,V_{1}, \ldots , V_{d}$ defined in the preceding paragraph,
$R \in \lt$ unknown). Of course, if $|q|< \sqrt{2} -1$, and if $\rho$ is the
solution of (1.4) given by \cite{JWS}, then $\fiq ( \rho )$
satisfies (1.5); in addition, $\fiq ( \rho )$ leaves
invariant each subspace of ${\cal T}$ spanned by tensors of a given
length (this follows immediately from the fact that $\rho$ can be
obtained by doing iterations in (1.4), starting with $\rho_{1} =p$).

We shall prove the following:

$\ $

{\bf Theorem} $1^{o}$ For every $-1<q<1$, there exists a unique
positive operator $R \in \lt$ which satisfies (1.5) and leaves
invariant each subspace of $\cal T$ spanned by tensors of a
given length.

$2^{o}$ For every $-1<q<1$, there exists a canonical unitary
$U : \tq \rightarrow {\cal T}$ which intertwines $R \in \lt$
defined above with
${\left( \sum_{i=1}^{d} A_{i}^{*} A_{i} \right)}^{1/2} \in \ltq$.
Moreover, we have that
\begin{equation}
UA_{i}U^{*} = V_{i}R, \ \ 1 \leq i \leq d.
\end{equation}

$3^{o}$ For every $-1<q<1$, the C*-algebra
$U \rqiu U^{*} \subset \lt$ contains $\ro$.

$4^{o}$ For $q$  satisfying:
\begin{equation}
q^{2} \ < \ 1-2|q|+2|q|^{4}-2|q|^{9}+ \cdots +2(-1)^{k}|q|^{k^{2}}+ \cdots
\end{equation}
the inclusion $U \rqiu U^{*} \subset \ro$
also holds, and hence $\rqiu$ is isomorphic to the extension by the
compacts of the Cuntz algebra.

$\ $

The inequality (1.7) gives
for $|q|$ a bound of around 0.44.
Calculations by computer indicate that Proposition 5.2 of the
paper, which has the last part of the Theorem as a corollary,
actually works (and gives $U \rqiu U^{*} = \ro$)
for $|q|$ up to a bound somewhere between 0.455 and 0.47.

The above theorem also gives some information on $\rqiu$ for larger
values of $q$ - for instance that
$\rqiu$ contains the compact operators on the
twisted Fock space $\tq$ for all the values of the parameter.
(We suspect this was known by
people working on the problem, altough it has not yet appeared in
writing.)

The paper is subdivided into sections as follows: in Section 2 we
review some basic facts about the twisted Fock space, and fix our
notations. In Section 3 we introduce the canonical unitary
$U: \tq \rightarrow {\cal T}$ and prove the first two assertions of
the above Theorem. In Section 4 we show that
$U \rqiu U^{*} \supseteq \ro$, and in Section 5 we prove the opposite
inclusion for $q$ satisfying (1.7).

$\ $

{\bf Acknowledgements:} We  would like to express our warmest thanks
to Victor Nistor for enjoyable and far-ranging conversations. We also
thank Gabriel Nagy for useful discussions at an early stage of this
work, and Roland Speicher for bringing \cite{Z} to our attention.

$\ $

$\ $

\section{Preliminaries}

$\ \ \ \ \ ${\bf 2.1 The Fock representation} If $\xi_{1}, \ldots , \xi_{d}$
is an orthonormal basis of the space $\cal H$ appearing in (1.1)
\setcounter{equation}{0}
(which is fixed throughout the paper, and has finite dimension
$d \geq 2$), then an orthonormal basis of the Fock space on
$\cal H$, ${\cal T} = {\bf C} \bigoplus \left(
\bigoplus_{n=1}^{\infty} {\cal H}^{\otimes n} \right)$, is:
\begin{equation}
\{ \Omega \} \bigcup \{ \sirxii \ | \ n \geq 1, \siri \};
\end{equation}
$\Omega$ in (2.1) (the vacuum vector) is 1 in the first
summand, $\bf C$, in the expression of $\cal T$. For every $n \geq 0$ we
shall denote by $\vn \subset {\cal T}$ the subspace spanned by
tensors of length $n$ in (2.1) ($\vz = {\bf C} \Omega$, by
convention; clearly dim $\vn = d^{n}, \ n \geq 0$, and
${\cal T} = \bigoplus_{n=0}^{\infty} \vn$, orthogonal direct sum). The
annihilation operators $V_{1}, \ldots , V_{d}$ involved
in equation (1.5) are determined by
\begin{equation}
V_{i} \Omega = 0, \ \  V_{i} ( \sirxii ) \ = \ \delta_{i,i_{1}}
\sirxiidoi ;
\end{equation}
their adjoints are the corresponding creation operators:
\[
V_{i}^{*} \Omega = \xi_{i}, \ \ V_{i}^{*} ( \sirxii ) \ = \
\xi_{i} \otimes \sirxii .
\]
$\ \ \ \ $Let us now pick a $q \in (-1,1)$. One defines recursively a
$q$-inner product $< \cdot \ | \ \cdot >_{q}$ on the subspaces
$\vn \subset {\cal T}, \ (n \geq 0)$, as follows: on $\vz$,
$< \cdot \ | \ \cdot >_{q}$ is determined by
$< \Omega \ | \ \Omega >_{q} = 1$;
then for every $n \geq 1$, one puts:
\[
< \sirxii \ | \ \sirxij >_{q} \ =
\]
\[
= \ \sum_{k=1}^{n} q^{k-1} \delta_{i_{1},j_{k}}
< \sirxiidoi \ | \ \sirxijhat >_{q}
\]
(where
$1 \leq i_{1}, \ldots , i_{n},j_{1}, \ldots ,j_{n} \leq d$, and
the hat on $\xi_{j_{k}}$ means that
$\xi_{j_{k}}$ is deleted from the tensor). We shall denote
$\vn$, considered with the $q$-inner product, by $\vnq$.
The natural basis of $\vn$, consisting of tensors of length
$n$ from (2.1), will be no longer orthogonal for
$< \cdot \ | \ \cdot >_{q}$ (unless $q$=0 or $n  \leq 1$); the
point is, however, that the Gramm matrix
$\Gamma_{n}$ of $q$-inner products
of elements from this basis remains positive and
non-degenerate (\cite{BS}, \cite{F}, \cite{Z}). Hence
one can define the Hilbert
space $\tq  =  \bigoplus_{n=0}^{\infty}  \vnq$
(orthogonal direct sum);
this is the twisted Fock space. The operators
$A_{1}, \ldots , A_{d}$ of (1.3) act on $\tq$ by:
\[
A_{i} \Omega \ = \ 0, \ \ A_{i} ( \sirxii ) \ = \
\sum_{k=1}^{n} q^{k-1} \delta _{i,i_{k}} \sirxiihat ;
\]
their adjoints are the corresponding creation operators:
\[
A_{i}^{*} \Omega \ = \ \xi_{i}, \ \ A_{i}^{*} ( \sirxii ) \ = \
\xi_{i} \otimes \sirxii .
\]
An important role in what follows will be played by the operator
\begin{equation}
M \ = \ \sum_{i=1}^{d} A_{i}^{*} A_{i} \in \ltq .
\end{equation}
Note that $M$ leaves invariant every subspace $\vnq \subset \tq$
spanned by tensors of length $n$. We shall denote by
$M_{n} \in  \lvnq$ the operator induced by $M$ on $\vnq$, and by
$[M_{n}] \in {\mbox{Mat}}_{d^{n}} ( {\bf C})$ the matrix of
$M_{n}$ with respect to the natural basis
$\{ \sirxii \ | \ \siri \}$ of $\vnq$ (ordered lexicographically,
for instance).  It is important to make distinction between
$M_{n}$ and $[M_{n}]$, since, due to the non-orthogonality of
the natural basis of $\vnq$, the matrix $[M_{n}]$ is generally
non-selfadjoint, altough $M_{n}$ itself is positive.

In general, if $X$ is in either of $\lvnq$, $\lvnqvn$,
$\lvnvnq$ ($n \geq 0, \ -1 < q < 1$), we shall denote by
$[X]$ its matrix with respect to the natural basis of its
domain and codomain. We have the usual rules
$[\alpha X + \beta Y ] = \alpha [X] + \beta [Y]$,
$[XY] = [X][Y]$, but computing $[X^{*}]$ needs a correction
with the Gramm matrix $\Gamma_{n}$ of $q$-inner products
of vectors from the natural basis of $\vnq$:
\begin{equation}
[X^{*}] \ = \ \left\{  \begin{array}{ll}
\Gamma_{n}^{-1} [X]^{*} \Gamma_{n} & \mbox{   if } X \in \lvnq  \\
\Gamma_{n}^{-1} [X]^{*}           & \mbox{   if } X \in \lvnqvn  \\
{[X]}^{*} \Gamma_{n}                & \mbox{   if } X \in \lvnvnq  \\
\end{array}  \right.
\end{equation}
(where $[X]^{*}$ is the conjugated-transpose of the matrix
$[X]$).

Though there were objects defined in this subsection (or in
the Introduction) which depend implicitly on the parameter $q$,
but do not have this dependence reflected in their notations,
we hope that this will not create any confusion in what follows.
The next list may also be of some help.

$\ $

\begin{tabular}{|c|c|}  \hline
Depend on $q$      & Don't depend on $q$     \\
			& (correspond to $q=0$)  \\   \hline
$a_{1}, \ldots , a_{d}$   &  $v_{1}, \ldots , v_{d}$       \\
$\tq = \bigoplus_{n=0}^{\infty} \vnq$   &
                        ${\cal T} = \bigoplus_{n=0}^{\infty} \vn$  \\
$A_{1}, \ldots , A_{d}$      &   $V_{1}, \ldots , V_{d}$        \\
$\bigoplus_{n=0}^{\infty} M_{n} =M= \sum_{i=1}^{d} A_{i}^{*} A_{i}$  &
		 $P= \sum_{i=1}^{d} V_{i}^{*} V_{i}$       \\
$\Gamma_{n} , \ n \geq 0$     &                            \\
$U= \bigoplus_{n=0}^{\infty} U_{n}$ (see Def. 3.2)   &       \\
$R= \bigoplus_{n=0}^{\infty} R_{n}$ (see Def. 3.3)   &      \\  \hline
\end{tabular}

$\ $

$\ $

{\bf 2.2 The actions of symmetric groups} Let $n \geq 1$ be an
integer, and let $S_{n}$ be the group of all permutations of
$\{ 1, \ldots ,n \}$. For $1 \leq k \leq l \leq n$,
we shall denote the cycle

$\left(  \begin{array}{ccccc}
k   & k+1 & \ldots & l-1 & l  \\
k+1 & k+2 &        & l   & k
\end{array}  \right) \in S_{n}$
by $\ciclukl$ (if $k=l$, then $\ciclukl$ is by
convention the unit of $S_{n}$). For every $-1 < q < 1$, we
have a natural representation by invertible operators
$\pi_{n,q} : S_{n} \rightarrow \lvnq$, determined by:
\begin{equation}
\pi_{n,q} (s) ( \sirxii ) \ = \
\xi_{i_{s^{-1} (1)}} \otimes \cdots \otimes
\xi_{i_{s^{-1} (n)}} .
\end{equation}
$\pi_{n,q}$ extends to a representation of ${\bf C}[S_{n}]
= \mbox{C*} (S_{n})$ on $\vnq$, still denoted by $\pi_{n,q}$.
This is  generally not a $\star$-representation; however, let us
point out that since invertible elements of C*($S_{n}$) must map
to invertible operators on $\vnq$, it is true that the spectrum
of $\pi_{n,q}(x)$ is contained in the one of $x$, for every $x$
in C*($S_{n}$).
This will be useful for studying the spectrum of $M$, since,
as shown
by a moment's reflection, we have:
\begin{equation}
M_{n} \ = \ \pi_{n,q} \left( \sum _{k=1}^{n}
q^{k-1} \cicluunuk \right) ,  \ \  n \geq 1 .
\end{equation}

$\ $

\section{The canonical unitary $U: \tq \rightarrow {\cal T}$}

$\ \ \ \ \ ${\bf 3.1 Lemma} For every $q \in (-1,1)$
\setcounter{equation}{0}
and $n \geq 1$ we have
\begin{equation}
\Gamma_{n} \ = \ \left(
\begin{array}{ccc}
\Gamma_{n-1}  &         &    0           \\
	      & \ddots  &               \\
  0	      &         &  \Gamma_{n-1} \\
\end{array}
\right)  [M_{n}]
\end{equation}
(equality in ${\mbox{Mat}}_{d^{n}}({\bf C})$, with $\Gamma_{n}$,
$\Gamma_{n-1}$, $[M_{n}]$ as in Section 2.1).

$\ $

{\bf Proof} Consider the representation
$\pi_{n,q} : \mbox{C*}(S_{n}) \rightarrow \lvnq$ defined in
Section 2.2. From Lemma 3 of \cite{BS} it follows that
$\pi_{n,q} \left( \sum_{s \in S_{n}} q^{inv(s)} s \right)$
has matrix $\Gamma_{n}$ in the natural basis of $\vnq$, where
$inv(s)$ is the number of inversions of $s \in S_{n}$.
Similarly, the matrix of
$\pi_{n,q} \left( \sum_{t \in S_{n},t(1)=1} q^{inv(t)} t \right)$ is
$\left( \begin{array}{ccc}
\Gamma_{n-1}  &         &      0         \\
	      & \ddots  &               \\
       0      &         &  \Gamma_{n-1} \\
\end{array} \right)$.
Taking also (2.6) into account, we see
that (3.1) is implied by:
\[
\sum_{s \in S_{n}} q^{inv(s)} s  \ = \
\left( \sum_{t \in S_{n},t(1)=1} q^{inv(t)} t \right)
\left( \sum _{k=1}^{n} q^{k-1} \cicluunuk \right);
\]
but this in turn comes out, exactly as in Proposition 1
of \cite{Z}, from the fact that every $s \in S_{n}$
can be uniquely decomposed as a product
$t \cicluunuk$ with $k \geq 1,$
$t \in S_{n}, \ t(1)=1$, and that in this decomposition we have
$inv(s)=inv(t)+(k-1)$. {\bf QED}

$\ $

{\bf 3.2 Proposition and Definition} Let $q$ be in (-1,1).
Let $U_{0} : \vzq \rightarrow \vz$ be defined by
$U_{0} \Omega = \Omega$, and then define recursively, for
$n \geq 1$:
\begin{equation}
U_{n} \ = \ (I \otimes U_{n-1})M_{n}^{1/2} : \vnq
\rightarrow \vn .
\end{equation}
(In the last formula, $I \otimes U_{n-1} \in \lvnqvn$
sends the tensor $\sirxii \in \vnq$ into
$\xi_{i_{1}} \otimes (U_{n-1} ( \sirxiidoi )) \in \vn$.)
Then $U_{n}$ is unitary, for every $n \geq 0$, and thus
$U = \bigoplus_{n=0}^{\infty} U_{n}$ is a unitary between
$\tq$ and $\cal T$.

$\ $

{\bf Proof} The fact that $U_{n}$ is unitary is equivalent to
\begin{equation}
[U_{n}]^{*} [U_{n}] \ = \ \Gamma_{n}
\end{equation}
(because $[U_{n}^{*} U_{n}] = [U_{n}^{*}][U_{n}]
\begin{array}[b]{c}   (2.4)   \\
		     =
\end{array}
\Gamma_{n}^{-1} [U_{n}]^{*} [U_{n}]$ ).
We prove (3.3) by induction on $n$. The case $n=0$ is clear.
The induction step ($n-1 \Rightarrow n$): obviously
\[
[U_{n}] \ = \ [I \otimes U_{n-1}][M_{n}^{1/2}] \ = \
\left( \begin{array}{ccc}
[U_{n-1}]     &         &               \\
	      & \ddots  &               \\
	      &         &  [U_{n-1}]  \\
\end{array} \right) [M_{n}^{1/2}],
\]
hence:
\[
[U_{n}]^{*}[U_{n}] \ = \ [M_{n}^{1/2}]^{*}
\left( \begin{array}{ccc}
[U_{n-1}]^{*}[U_{n-1}]  &         &               \\
	      & \ddots  &               \\
	      &         &  [U_{n-1}]^{*}[U_{n-1}]  \\
\end{array} \right) [M_{n}^{1/2}]
\]
\[
= \ [M_{n}^{1/2}]^{*}
\left( \begin{array}{ccc}
\Gamma_{n-1}  &         &               \\
	      & \ddots  &               \\
	      &         &  \Gamma_{n-1} \\
\end{array} \right) [M_{n}^{1/2}];
\]
replacing $[M_{n}^{1/2}]^{*}$ from relation (2.4) and using
Lemma 3.1 we get that this equals:
\[
( \Gamma_{n} [M_{n}^{1/2}] \Gamma_{n}^{-1} )
( \Gamma_{n} [M_{n}^{-1}] ) [M_{n}^{1/2}]
= \ \Gamma_{n} [M_{n}^{1/2} M_{n}^{-1} M_{n}^{1/2}] \ = \
\Gamma_{n} . \ \
{\bf QED}
\]

$\ $

{\bf 3.3 Remark} The construction of the above unitary $U$ would still
work if $d$ (the dimension of the separable Hilbert space $\cal H$
of (1.1)) were infinite.
In this case, the spaces $(\vn)_{n=0}^\infty$ and $(\vnq)_{n=0}^\infty$
we are working with would of course no longer be finite dimensional.
However, for every sequence $\alpha=(\alpha_1,\alpha_2,\alpha_3,\ldots)$
of non--negative integers with $|\alpha|=\sum_{i=1}^\infty\alpha_i<\infty$,
let us denote by $\va$ (respectively $\vaq$) the subspace of $\cal T$
(respectively $\tq$) generated by the tensors
$\xi_{i_1}\otimes\cdots\otimes\xi_{i_n}$ having the property that among
$i_1,\ldots,i_n$ there are $\alpha_1$ of $1$, $\alpha_2$ of $2$,
$\alpha_3$ of 3, $\ldots$.
Then for every $n\ge0$ we have $\vn=\bigoplus_{|\alpha|=n}\va$,
$\vnq=\bigoplus_{|\alpha|=n}\vaq$, orthogonal direct sums, and each
space $\va$ and $\vaq$ is finite dimensional even when $d$
is infinite.
Moreover, one can rewrite the construction of the above
$U:\tq\rightarrow{\cal T}$ by defining a family of unitaries
$(U_\alpha:\vaq\rightarrow\va)_\alpha$, by induction on $|\alpha|$, and
then putting $U=\bigoplus_\alpha U_\alpha$.
(When doing so, $M_n$ of formula (3.2) should be replaced by
$M_\alpha\in\lvaq$, determined by
$$ M_\alpha(\xi_{i_1}\otimes\cdots\otimes\xi_{i_n})
  = \sum_{k=1}^nq^{k-1}\xi_{i_k}\otimes\xi_{i_1}\otimes\cdots\otimes
    \xi_{i_{k-1}}\otimes\xi_{i_{k+1}}\otimes\cdots\otimes\xi_{i_n}, $$
for $\xi_{i_1}\otimes\cdots\otimes\xi_{i_n}\in\vaq$.)
This way of proving the preceding proposition still works for
$d=\infty$, (but in what follows, we are only concerned with the case
of finite $d$).

$\ $

{\bf 3.4 Definition and Proposition} Let $q$ be in (-1,1),
and define:
\begin{equation}
R \ = \ U M^{1/2} U^{*} \in \lt .
\end{equation}
Then: $1^{o}$ We have
\begin{equation}
U A_{i} U^{*} \ = \ V_{i} R, \ \ 1 \leq i \leq d,
\end{equation}
and as a consequence $V_{1}R, \ldots , V_{d}R$ satisfy the
\qcom.

$2^{o}$ $R$ is the unique positive operator on the (non-deformed)
Fock space $\cal T$ which satisfies equation (1.5), and which
leaves invariant each subspace of $\cal T$ spanned by tensors of a
given length.

$\ $

{\bf Proof} $1^{o}$ We pick $1 \leq i \leq d$ and check the
equality $A_{i}^{*} = M^{1/2} U^{-1} V_{i}^{*} U$ (obviously
equivalent to (3.5)) on vectors of the natural basis of $\tq$.
Note first that
$A_{i}^{*} \Omega = \xi_{i} = M^{1/2} U^{-1} V_{i}^{*} U \Omega$
(the second equality following from the facts, easy to check, that
$U \xi_{i} = \xi_{i}$ and that $M_{1}$ is the identity on the space
${\cal V}_{1,q}$). Next, for a tensor $\sirxii \in \tq$ we have:
\[
M^{1/2} U^{-1} V_{i}^{*} U \sirxii \ = \
M_{n+1}^{1/2} U_{n+1}^{-1} V_{i}^{*} U_{n} \sirxii .
\]
But $M_{n+1}^{1/2} U_{n+1}^{-1} = (I \otimes U_{n})^{-1}$
(by (3.2)), which implies that the coincidence of $A_{i}^{*}$
and $M^{1/2} U^{-1} V_{i}^{*} U$ on $\sirxii$ is equivalent to:
\[
(I \otimes U_{n} ) A_{i}^{*} \sirxii \ = \
V_{i}^{*} U_{n} \sirxii .
\]
The last equality is, however, obvious from the definitions
of $A_{i}^{*}$ and $V_{i}^{*}$.

$2^{o}$ It is clear that $R$ is positive and leaves invariant
every $\vn$ (= subspace spanned by tensors of length $n$).
The fact that $R$ satisfies equation (1.5) follows from (3.5)
via exactly the argument preceding relation (10) of \cite{JWS}.

If $\widetilde{R} \in \lt$ is an operator sharing the above properties,
then (as it immediately comes out of (1.5)), $\widetilde{R} \Omega =0$
and, for $m,n \geq 1$,
$1 \leq i_{1}, \ldots , i_{m},j_{1}, \ldots ,j_{n}  \leq d$:
\[
< {\widetilde{R}}^{2} ( \xi_{i_{1}} \otimes \cdots \otimes
\xi_{i_{m}} ) \ | \ \sirxij > \ =
\]
\[
= \ \delta_{m,n} \left( \delta_{i_{1},j_{1}} \cdots
\delta_{i_{n},j_{n}} + q < V_{j_{1}} \widetilde{R}
\xi_{i_{2}} \otimes \cdots \otimes \xi_{i_{m}} \ | \
V_{i_{1}} \widetilde{R}
\xi_{j_{2}} \otimes \cdots \otimes \xi_{j_{n}} > \right) .
\]
 From the last equation it is clear that, for every $n \geq 1$,
$\widetilde{R} | {\cal V}_{n-1}$ determines ${\widetilde{R}}^{2} | \vn$
(and hence $\widetilde{R} | \vn$ too, by taking a square root).
Thus an induction argument shows that $\widetilde{R} | \vn = R | \vn$
for every $n \geq 0$, and we get $\widetilde{R} = R$.  {\bf QED}

$\ $

$\ $

\section{$U \rqiu U^{*} \supseteq \ro$}

$\ \ \ \ \ ${\bf 4.1 Lemma} Let $q$ be in (-1,1),
let $R$ be as in Definition 3.4,
\setcounter{equation}{0}
and denote, for every $n\geq 0$, $R| \vn$ by $R_{n}$. Then we have
\begin{equation}
\left( \frac{1}{1-|q|} \prinf \right) I \leq R_{n}^{2} \leq
\frac{1}{1-|q|} I
\end{equation}
(inequality in $\lvn$), for every $n \geq 1$.

$\ $

{\bf Proof} Since $R_{n}^{2}$ is conjugate
to $M_{n}$ by $U_{n}$, it suffices to prove the analogue in
$\lvnq$ of (4.1), with $R_{n}^{2}$ replaced by $M_{n}$.
This comes, clearly, to showing that the spectral radii of
$M_{n}$ and $M_{n}^{-1}$ are not greater than $1/(1-|q|)$ and
$(1-|q|) \prinfinv$, respectively. Recalling the considerations
of Section 2.2, and dominating the spectral radius of an element
$x \in \mbox{C*}(S_{n})$ by $||x||$, we see that it will suffice
to prove:
\begin{equation}
\left\{  \begin{array}{l}
|| \xn || \leq \frac{1}{1-|q|}    \\
|| {\left( \xn \right)}^{-1} || \leq (1-|q|) \prinfinv   \\
\end{array}  \right.  .
\end{equation}
The first inequality in (4.2) is obvious. The proof of the
second one (which must, of course, contain a proof of the
invertibility of $\xn$), is obtained from a decomposition into
factors performed in the spirit of \cite{Z}, Proposition 2.
More precisely, one notes first the commutation relation:
\[
(1 \rightarrow j)(1 \rightarrow k) \ = \
(2 \rightarrow k)(1 \rightarrow j-1), \ \ 2 \leq j \leq k \leq n,
\]
which immediately implies the identity:
\[
\left( \sum_{k=1}^{m} q^{k-1} (1 \rightarrow k) \right)
\left( I-q^{m-1} (1 \rightarrow m) \right) \ = \
\]
\begin{equation}
= \ \left( I- q^{m} (2 \rightarrow m) \right)
\left( \sum_{k=1}^{m-1} q^{k-1} (1 \rightarrow k) \right) ,
2 \leq m \leq n.
\end{equation}
Multiplying (4.3) by
${\left( I-q^{m-1} (1 \rightarrow m) \right)}^{-1}$ on the right,
and using induction on $m$ ($1 \leq m \leq n$), one obtains:
\begin{equation}
\xn \ = \ \prod_{j=0}^{n-2}
\left( I-q^{n-j} (2 \rightarrow n-j) \right)
\prod_{j=1}^{n-1}
{\left( I-q^{j} (1 \rightarrow j+1) \right)}^{-1} .
\end{equation}
Taking the inverse in (4.4), and after that taking norms and
doing straightforward majorizations, we get the second inequality
(4.2) (the norm of
${\left( I-q^{n-j} (2 \rightarrow n-j) \right)}^{-1} =
\sum_{k=0}^{\infty} q^{(n-j)k} (2 \rightarrow n-j)^{k}$
is dominated by $1/(1-|q|^{n-j})$).  {\bf QED}

$\ $

{\bf 4.2 Remark} $\prinf$ appearing in (4.1) is strictly
positive (as it is well-known, $\sum_{k=1}^{\infty} |q|^{k} < \infty$
implies $\prod_{k=1}^{\infty} (1+|q|^{k}) < \infty$ and
$\prod_{k=1}^{\infty} (1-|q|^{k}) > 0$). We take this occasion to
mention the following identity (due to Gauss - see Corollary 2.10
of \cite{Andrews}):
\begin{equation}
\prod_{k=1}^{\infty} \frac{1-q^{k}}{1+q^{k}} \ = \
\sum_{k=- \infty}^{\infty} (-1)^{k} q^{k^{2}} ,  \ \ 0 \leq q < 1.
\end{equation}

$\ $

{\bf 4.3 The proof} of the inclusion
contained in the title of this section
is now immediate. Indeed, what we need to show is that
$V_{1}, \ldots , V_{d} \in U \rqiu U^{*}$. We know, from the
relations (3.5), that $V_{1}R, \ldots , V_{d}R \in U \rqiu U^{*}$,
where $R= UM^{1/2} U^{*}$ also belongs to $U \rqiu U^{*}$. But
now, $R$ splits as $\bigoplus_{n=0}^{\infty} R_{n}$, with
$R_{0}=0$ on ${\bf C} \Omega$ and $R_{n} \in \lvn$ satisfying
``the square root'' of inequality (4.1) for $n \geq 1$; this
immediately implies that Ker$R = {\bf C} \Omega$, and that 0
is isolated in the spectrum of $R$. Denoting by $\phi$ the
(continuous) function on the spectrum of $R$ which sends 0
into 0 and $\alpha \neq 0$ into $1/ \alpha$, we then have that
$R \phi (R)$ is the projection onto
${\cal T} \ominus {\bf C} \Omega$, and hence that indeed
$V_{i}=(V_{i}R) \phi (R) \in U \rqiu U^{*}, \ 1 \leq i \leq d$.
{\bf QED}

$\ $

$\ $

\section{$U \rqiu U^{*} \subseteq \ro$ for $q$ satisfying (1.7)}

$\ \ \ \ \ $Let $q$ be in (-1,1), and let $R \in \lt$ be as in Definition 3.4.
\setcounter{equation}{0}
Due to the relations (3.4) and (3.5), it is obvious that the inclusion
``$U \rqiu U^{*} \subseteq \ro$'' is equivalent to ``$R \in \ro$''.

Denote now, as in Lemma 4.1, by $R_{n}$ the operator induced by
$R$ on the invariant subspace $\vn \subset {\cal T}, \ (n \geq 0)$,
and define, for every $n \geq 1$:
\begin{equation}
X_{n} \ = \ R_{0} \oplus R_{1} \oplus \cdots R_{n} \oplus
(R_{n} \otimes I) \oplus (R_{n} \otimes I \otimes I) \oplus \cdots
\in \lt .
\end{equation}

It is immediate (from Definitions 3.2 and 3.4) that $R_{0}=0$ and
$R_{1}$ is the identity operator on ${\cal V}_{1}$. Moreover, for
$n \geq 1$, restricting the two sides of the equation (1.5) to
the subspace $\vnp$ gives:
\begin{equation}
R_{n+1}^{2} \ = \ I_{n+1} + q \sum_{i,j=1}^{d}
V_{j}^{*} R_{n} V_{i}^{*} V_{j} R_{n} V_{i}
\end{equation}
(where $I_{n+1}$ is the identity operator on $\vnp$, and we view
$V_{i} \in {\cal L} ( {\cal V}_{n+1}, {\cal V}_{n} )$,
$V_{j} \in {\cal L} ( {\cal V}_{n}, {\cal V}_{n-1} )$,
$V_{i}^{*} \in {\cal L} ( {\cal V}_{n-1}, {\cal V}_{n} )$,
$V_{j}^{*} \in {\cal L} ( {\cal V}_{n}, {\cal V}_{n+1} )$ ).
Using all these facts, it is easy to check that the $X_{n}$'s
defined in (5.1) satisfy:
\begin{equation}
X_{n+1}^{2} \ = \ P+q \sum_{i,j=1}^{d} V_{j}^{*}X_{n}V_{i}^{*}
V_{j}X_{n}V_{i}, \ \ n \geq 1.
\end{equation}
Thus $(X_{n})_{n=1}^{\infty}$ is the Fock representation of a
sequence of iterates as considered in \cite{JWS}, which begins
with $X_{1}=P$. Note that in this particular situation, the
iterates can be defined with no restriction on $q \in (-1,1)$.

$\ $

{\bf 5.1 Lemma} Let $q$ be in (-1,1), let
$R = \bigoplus_{n=0}^{\infty} R_{n}$ be as above and denote,
for every $n \geq 1$, by $\alpha_{n} (q)$ the smallest eigenvalue
of $R_{n}^{2}$ ($\alpha_{n} (q) > 0$
by Lemma 4.1). Then  the
$(X_{n})_{n=1}^{\infty}$ defined in (5.1) satisfy:
\begin{equation}
||X_{n+2} - X_{n+1}|| \ \leq \
\frac{|q|}{\sqrt{(1-|q|) \mbox{min} ( \alpha_{n+1} (q), \alpha_{n+2} (q) ) }}
||X_{n+1} - X_{n}||, \ \ n \geq 1 .
\end{equation}

$ \ $

{\bf Proof} The argument will consist in combining Lemma 8 of
\cite{JWS} with the particular form given to the iterates in (5.1).
It is immediate (from (5.1)) that
$|| X_{n} - X_{n+1}||=||(R_{n} \otimes I) - R_{n+1}||$, and we
shall examine the latter quantity.

Observe that because of the obvious relations
$R_{n}V_{i}=V_{i}(R_{n} \otimes I)$,
$V_{j}^{*} R_{n} = (R_{n} \otimes I)V_{j}^{*}$,
the equation (5.2) can be rewritten
\begin{equation}
R_{n+1}^{2} \ = \ I_{n+1} + q (I \otimes R_{n}) T_{n+1}
(I \otimes R_{n}),
\end{equation}
where $T_{n+1}$ is the operator induced on $\vnp$ by
$T= \sum_{i,j=1}^{d} V_{j}^{*}V_{i}^{*}V_{j}V_{i}$.
Since clearly $T_{n+1} \otimes I = T_{n+2}$, the last
equality gives, when tensored with $I$ on the right:
\begin{equation}
R_{n+1}^{2} \otimes I \ = \ I_{n+2} +
q (I \otimes R_{n} \otimes I) T_{n+2}
(I \otimes R_{n} \otimes I).
\end{equation}
Thus, if in $R_{n+2}^{2}-(R_{n+1}^{2} \otimes I)$ we
replace $R_{n+2}^{2}$ using the analogue of (5.5) (for $n+2$) and
$R_{n+1}^{2} \otimes I$ using (5.6), we obtain, by taking
norms:
\[
||R_{n+2}^{2}-(R_{n+1}^{2} \otimes I)|| \ =
\]
\[
= \ |q| \ ||(I \otimes R_{n+1}) T_{n+2}
(I \otimes R_{n+1} - I \otimes R_{n} \otimes I) +
(I \otimes R_{n+1} - I \otimes R_{n} \otimes I)
T_{n+2} (I \otimes R_{n} \otimes I)||
\]
\[
\leq \ |q| (||R_{n}||+||R_{n+1}||)
||R_{n+1}-(R_{n} \otimes I)||.
\]

On the other hand, $R_{n+2}^{2} \geq  \alpha_{n+2} (q) I_{n+2}$,
$R_{n+1}^{2} \otimes I \geq \alpha_{n+1} (q) I_{n+2}$,
hence Lemma 8 of \cite{JWS} gives that
\[
||R_{n+2}-(R_{n+1} \otimes I)|| \ \leq \
\frac{1}{2 \sqrt{ \mbox{min} ( \alpha_{n+1} (q), \alpha_{n+2} (q) ) }}
||R_{n+2}^{2}-(R_{n+1}^{2} \otimes I)||.
\]
Combining this with the bound obtained for
$||R_{n+2}^{2}-(R_{n+1}^{2} \otimes I)||$ (in which $||R_{n}||$,
$||R_{n+1}||$ are majorized, by Lemma 4.1,
with $1/ \sqrt{1-|q|}$ ), we get (5.4).
{\bf QED}

$\ $

$\ $

{\bf 5.2 Proposition} Let $q$ be in (-1,1), and let
$R= \bigoplus_{n=0}^{\infty} R_{n}$ and
$( \alpha_{n} (q))_{n=1}^{\infty}$ be as in
Lemma 5.1. If
\begin{array}[t]{c}   \mbox{lim inf}  \\
		       n \rightarrow \infty \\
\end{array}
$ \ \alpha_{n} (q) \ > \ q^{2}/(1-|q|)$, then
$R \in \ro$
(and hence $\rqiu$ is isomorphic to the extension by the
compacts of the Cuntz algebra).

$\ $

{\bf Proof} Consider the sequence
$(X_{n})_{n=1}^{\infty}$ defined in (5.1). From Lemma 5.1
and the ratio test it follows that
$\sum_{n=1}^{\infty} ||X_{n+1}-X_{n}|| < \infty$,
hence this sequence is norm convergent. Each $X_{n}$
is in $\ro$ (as it is clear by using (5.3) and an induction
argument), hence the limit is in $\ro$, too. But the limit
can only be $R$ (indeed, it is obvious from (5.1) that
$X_{n}$ converges to $R$ in the strong operator topology). {\bf QED}

$\ $

$\ $

{\bf 5.3 Corollary} If $q$ satisfies (1.7), then $\rqiu$ is
isomorphic to the extension by the compacts of the Cuntz algebra.

$\ $

{\bf Proof} Since, by Lemma 4.1,
\[
\begin{array}[t]{c}   \mbox{lim inf}  \\
		       n \rightarrow \infty \\
\end{array}
\ \alpha_{n} (q) \ \geq \ \frac{1}{1-|q|}
\prod_{k=1}^{\infty} \frac{1-|q|^{k}}{1+|q|^{k}} \ = \
\frac{1}{1-|q|}
\sum_{k=- \infty}^{\infty} (-1)^{k} |q|^{k^{2}} ,
\]
the last proposition can be applied to every $q$ satisfying (1.7).
{\bf QED}

$\ $

$\ $

{\bf 5.4 Remark} Truncating the series on the right-hand
side of (1.7) to its first two terms leads to the
bound $|q| < \sqrt{2} - 1$ found in \cite{JWS}.
The only positive root of the equation
$q^{2}=1-2q+2q^{4}-2q^{9}$ is 0.44005651..., hence
taking four terms of the series makes us sure that (1.7)
is fulfilled for $|q| \leq 0.44$; the improvement obtained
by considering further terms of the series is only in the
sixth significant figure of the numerical value of the bound.
Of course, further improvements can be obtained by giving
better estimates for
\begin{array}[t]{c}   \mbox{lim inf}  \\
		       n \rightarrow \infty \\
\end{array}
$\ \alpha_{n} (q)$.
Computer aided calculations of this quantity
indicate that the hypothesis of
Proposition 5.2 is still fulfilled for $|q|$=0.455; however,
it appears that a new idea is certainly needed in order to
reach, say, 0.47.

$\ $

$\ $

\pagebreak


\begin{thebibliography}{99}

\bibitem{Andrews} G.E. Andrews. The theory of partitions,
Addison-Wesley, 1976.

\bibitem{BS} M. Bo\.{z}ejko, R. Speicher.
An example of a generalized
Brownian motion,  part I in Commun. Math. Phys., vol.137(1991),
519-531, part II in Quantum Probability and Related Topics VII,
Proceedings New Delhi 1990, 67-77.

\bibitem{C} J. Cuntz. Simple C*-algebras generated by isometries,
Commun. Math. Phys., vol.57(1977), 173-185.

\bibitem{E} D.E. Evans. On ${\cal O}_{n}$,
Publ. RIMS Kyoto Univ., vol.16(1980), 915-925.

\bibitem{F} D. Fivel. Interpolation between Fermi and Bose
statistics using generalized commutators, Phys. Rev. Lett.,
vol.65(1990), 3361-3364.

\bibitem{G} O.W. Greenberg. Example of infinite
statistics, Phys. Rev. Lett., vol.64(1990), 705-708.

\bibitem{JWS} P.E.T. Jorgensen, L.M. Schmitt, R.F. Werner.
$q$-canonical commutation relations and stability of
the Cuntz algebra, to appear in the Pacific Journal of
Mathematics.

\bibitem{Z}  D. Zagier. Realizability of a model in
infinite statistics, Commun. Math. Phys., vol.147(1992), 199-210.
\end{thebibliography}
\end{document}